\documentclass[reqno, 11pt]{amsart}
\usepackage[utf8]{inputenc}
\usepackage[T1]{fontenc}
\usepackage{lmodern}
\usepackage{epsfig,amsmath, amsthm ,amsfonts,latexsym,amssymb}
\usepackage[abs]{overpic}
\usepackage[usenames,dvipsnames]{xcolor}
\usepackage{cite}
\usepackage[linktocpage=true]{hyperref}
\usepackage{caption}
\usepackage{subcaption}
\usepackage{dsfont}
\usepackage{mathtools}
\usepackage{tikz-cd}
\usepackage{comment}
\usepackage{enumitem}
\usepackage{marginnote}

\usepackage[all,cmtip]{xy} 

\usepackage{braket}
\usepackage{enumitem}
\usepackage{tikz}
\usetikzlibrary{shapes,arrows}
\usetikzlibrary{cd, positioning}
\usepackage{pgfplots}
\usepackage{adjustbox}

\usepackage{amssymb}
\usepackage{mathrsfs}

\usepackage{mathtools}
\usepackage{cite}
\usepackage{tikz-cd}

\hypersetup{
  colorlinks,
  citecolor=teal,
  linkcolor=purple,
  urlcolor=teal}
\theoremstyle{plain}

\newtheorem{theorem}{Theorem}

\newtheorem{dummy}{anything}[section]

\newtheorem{lemma}[dummy]{Lemma}

\newtheorem{proposition}{Proposition}[section]

\newtheorem{corollary}[dummy]{Corollary}

\theoremstyle{definition}
\newtheorem{definition}[dummy]{Definition}
\newtheorem*{defn}{Definition}
\newtheorem{example}[dummy]{Example}

\newtheorem{remark}[dummy]{Remark}
\newtheorem*{rem}{Remark}
\theoremstyle{remark}
\newcommand{\R}{\mathbb{R}}
\newcommand{\C}{\mathbb{C}}

    \usepackage{extpfeil}
    \usepackage[new]{old-arrows} 
    \usetikzlibrary{cd}
    \usetikzlibrary{calc}
    \usetikzlibrary{decorations.pathmorphing}
    \tikzset{curve/.style={settings={#1},to path={(\tikztostart)
        .. controls ($(\tikztostart)!\pv{pos}!(\tikztotarget)!\pv{height}!270:(\tikztotarget)$)
        and ($(\tikztostart)!1-\pv{pos}!(\tikztotarget)!\pv{height}!270:(\tikztotarget)$)
        .. (\tikztotarget)\tikztonodes}},
        settings/.code={\tikzset{quiver/.cd,#1}
            \def\pv##1{\pgfkeysvalueof{/tikz/quiver/##1}}},
        quiver/.cd,pos/.initial=0.35,height/.initial=0}

\def\R{\mathbb{R}}

\def\C{\mathbb{C}}

\DeclareFontFamily{U}{mathx}{}
\DeclareFontShape{U}{mathx}{m}{n}{<-> mathx10}{}
\DeclareSymbolFont{mathx}{U}{mathx}{m}{n}
\DeclareMathAccent{\widecheck}{0}{mathx}{"71}

\begin{document}

\title[A Poincar\'e--Birkhoff theorem for $C^0$-Hamiltonian maps]{A Poincar\'e--Birkhoff theorem for $C^0$-Hamiltonian maps}

\author{Arthur Limoge}

\address[A.\ Limoge]{Institut f\"ur Mathematik \\ Universit\"at Heidelberg \\ Germany}

\email{\href{mailto:arthur.limoge@gmail.com}{arthur.limoge@gmail.com}}

\author{Agustin Moreno}

\address[A.\ Moreno]{Institut f\"ur Mathematik \\ Universit\"at Heidelberg \\ Germany}

\email{\href{mailto:agustin.moreno2191@gmail.com}{agustin.moreno2191@gmail.com}}

\date{}





\maketitle

\begin{abstract}  We prove a higher-dimensional version of the well-known Poincaré--Birkhoff theorem, using Floer homology. We also prove a relative version for Lagrangian submanifolds. The motivation is finding periodic orbits and Hamiltonian chords in the circular, restricted three-body problem.
\end{abstract}

\tableofcontents

\section{Introduction}

The problem of finding closed periodic orbits in the planar, circular, restricted three-body problem (or planar CR3BP) goes back to ground-breaking work in celestial mechanics of Poincar\'e \cite{P87,P12}, building on work of G.W.\ Hill on the lunar problem \cite{H77,H78}. The basic scheme for his approach may be reduced to:
\begin{itemize}
    \item[(1)] Finding a global surface of section for the dynamics (thus reducing the problem to studying the first return map);
    \item[(2)] Proving a fixed point theorem for this return map.
\end{itemize}
This is the setting for the celebrated Poincar\'e--Birkhoff theorem, also known as Poincar\'e's last geometric theorem, proposed and confirmed in special cases by Poincar\'e and later proved in full generality by Birkhoff in \cite{Bi13}. It states that an area-preserving homeomorphism of the annulus that satisfies a \emph{twist} condition at the boundary (i.e.\ it rotates the two boundary components in opposite directions) admits at least two fixed points, and moreover admits infinitely many periodic points of arbitrary large period, as was later observed. 

The Poincar\'e--Birkhoff theorem was further refined by Arnold, who conjectured a higher-dimensional version \cite{Ar65,Ar78}, stating that a non-degenerate Hamiltonian on a closed symplectic manifold has at least as many $1$-periodic orbits as the sum of the Betti numbers of the manifold (the fixed point version of the Poincar\'e--Birkhoff theorem follows by a simple doubling trick). This result was proved by Floer \cite{F89} in the aspherical case, famously introducing the homology theory that bears his name, and went on to become a crucial cornerstone in modern symplectic geometry (many others later contributed to proving the Arnold conjecture in generality, which required the introduction of far-reaching foundational methods; see \cite{BX22} for the state of the art). The case for periodic points was conjectured by Conley \cite{Co} and proved by Ginzburg \cite{G10}, using Floer homological tools.

While Poincar\'e focused on the planar case of the CR3BP, where the massless particle is constrained to the plane, we drop this assumption and work with the \emph{spatial} problem, i.e.\ the particle moves in three-space. This problem is higher-dimensional, and therefore new tools need to be deployed. In order to address the existence of periodic orbits, we shall follow the same two-step scheme as Poincar\'e. 

The first step was successfully addressed by Moreno--van-Koert, who found an adapted open book in \cite{MvK2}, yielding a family of global hypersurfaces of section. As a result, they obtained that the (subcritical, near-primary) dynamics of the spatial CR3BP could be encoded by a \emph{global} Poincar\'e return map on a $4$-dimensional degenerate Liouville domain. It is worth emphasising that this result is fully non-perturbative, i.e.\ it holds for all the relevant range of parameters (as opposed to the analogous results in the planar CR3BP; cf.\ \cite{HSW19}). The second step was also addressed by Moreno--van-Koert, who first showed in \cite{MvK2} that this Poincar\'e return map is Hamiltonian in the interior, leading to their formulation of a generalised version of the Poincar\'e--Birkhoff theorem in \cite{MvK}. While this result made considerable headway by its generality, its applicability remained limited. One of the reasons why it has proved so hard to apply it to the CR3BP is that the twist condition as originally defined by Moreno--van-Koert is too restrictive (e.g.\ it is not open), and moreover does not take into account the boundary degeneracy of the symplectic form. In this paper, we address these two issues. 

Heuristically, the goal of this paper is to generalise the model of Moreno \& van Koert by getting rid of most of the technical assumptions necessary for their results to hold, thus closing a long and fruitful circle of ideas. Concretely, we shall address the second step in the above scheme, by proving a fixed-point theorem for suitable Hamiltonian maps on Liouville domains (which we call \emph{$C^0$-Hamiltonian twist maps}), i.e.\ a higher-dimensional generalisation of the classical Poincar\'e--Birkhoff theorem. This allows us to significantly relax the twist condition (which is now an \emph{open} condition), as well as address the boundary degeneracy. We emphasise that this is now the correct setup in which to describe the CR3BP from a Floer-homological point of view, since our Liouville domains of interest are degenerate.


We will also prove a relative version of the fixed-point theorem for Lagrangian submanifolds, which is a generalised version of the result in \cite{LM2}, with similar improvements. The motivation for this result is to obtain infinitely many spatial trajectories of the CR3BP which start and end orthogonally to the $xz$-plane – one of the symmetry planes of the problem –, an example of which being halo orbits, as well as infinitely many consecutive collision orbits (which normally can be perturbed to be near-collision orbits, useful for gravitational assist manoeuvres and low-altitude flybys in the context of trajectory design). Similarly, one hopes to obtain trajectories between different Lagrangians, e.g.\ starting normal to the $xz$-plane and ending in collision. This is related to the heuristic arguments outlined by the authors in \cite{LM, LM2}.

\smallskip

Intuitively, $C^0$-Hamiltonian twist maps are smooth and Hamiltonian in the interior of the Liouville domain on which they are defined, but only admit continuous extensions to the boundary, in such a way that there exists a generating Hamiltonian which is ``infinitely wrapping'' at the boundary. Alternatively, they may be described as maps on \emph{degenerate} Liouville domains (similar to a Liouville domain but where the symplectic form degenerates at the boundary), which are everywhere smooth and satisfy a weakened version of the twist condition as defined in \cite{MvK}; see Definition \ref{def:C0Hamtwistmaps}. These maps naturally arise as return maps on the page of an open book adapted to the dynamics (as in the CR3BP), and also in the context of billiards. We indeed expect these methods to work for billiard maps, although they \textit{a priori} yield weaker results than what can be proved with less technically demanding methods (e.g.\ Morse theory). We remark that the general fixed-point theorems that we derive are of a non-perturbative nature (similar to the original Poincar\'e--Birkhoff theorem). We will leave the problem of applying these results to the CR3BP to future work.

\medskip

\textbf{Fixed-point theorems.} Before stating the main theorems, we give a definition of $C^0$-Hamiltonian twist maps. Recall that a Liouville domain is an exact symplectic manifold with contact-type boundary. 

\begin{definition}\label{def:C0Hamtwistmaps} \textbf{($C^0$-Hamiltonian twist maps on Liouville domains)} Let $f: (W,\omega)\rightarrow (W,\omega)$ be a map on a Liouville domain, and let $\alpha$ be the contact form along $B=\partial W$. We say that $f$ is a \emph{$C^0$-Hamiltonian twist map} if the following conditions hold: 
\begin{itemize}
    \item\textbf{(Hamiltonian)} $f\vert_{\mbox{int}(W)}=\phi_H^1$ is generated by a $C^1$ Hamiltonian $H_t: \mbox{int}(W)\rightarrow \mathbb R$;
    \item\textbf{(Extension)} Both $f$ and the Hamiltonian $H_t$ admit $C^0$ extensions to the boundary, but not necessarily $C^1$ extensions; and
    \item\textbf{(Weakened Twist Condition)} Near the boundary $B$, the generating Hamiltonian vector field satisfies $h_t:=\alpha(X_{H_t})>0$, and $h_t\rightarrow +\infty$ as we approach $B$.
\end{itemize}

We say that the isotopy $H_t$ is \emph{infinitely strictly wrapping} or \emph{infinitely positively wrapping}, or simply \emph{infinitely wrapping}. 

\end{definition}

\textbf{Closed-string version.} We will then prove the following generalised version of the main theorem in \cite{MvK}.

\begin{theorem}[\textbf{Closed strings: Long interior orbits}]\label{thm:main_thm}
Let $f:(W,\omega)\rightarrow (W,\omega)$ be a $C^0$-Hamiltonian twist map on a Liouville domain. Assume the following:
\begin{itemize}
\item\textbf{(fixed points)} All fixed points of $f$ are isolated (i.e.\ finitely many);
\item \textbf{(First Chern class)} $c_1(W)=0$ if $\dim W \geq 4$;
\item\textbf{(Symplectic cohomology)} $SH^\bullet(W)$ is non-zero in infinitely many degrees.
\end{itemize}
Then $f$ has simple interior periodic points of arbitrarily large (integer, minimal) period.
\end{theorem}

\textbf{Improvements and comparison.} We will follow a similar strategy as in \cite{MvK}, which is adapted from the strategy in \cite{G10}, but we will use action growth as opposed to index growth to rule out undesired orbits. This means that we obtain the following improvements when compared to \cite{MvK}:

\begin{itemize} 
    \item We can relax the twist condition, which is now much weaker and an open condition (indeed, we only ask that the Hamiltonian vector field should be positively transverse to the contact structure at the boundary, as opposed to being parallel to the Reeb vector field). 
    \item We can relax the smoothness of the generating Hamiltonians from $C^2$ to $C^1$ (which is the minimal regularity possible in order to define a continuous Hamiltonian vector field), as we do not need to consider the linearisation of the Hamiltonian flow;
    \item We can get rid of the assumptions on index-positivity; 
    \item We can drop the assumption that the contact structure at the boundary should be globally trivial (needed to measure the $CZ$-index for index-positivity).
\end{itemize}

As a mild trade-off, we have to assume that symplectic homology is supported in infinitely many degrees, as opposed to simply being infinite-dimensional. In \cite{MvK}, this fact came for free under the index-positive assumption (here, note that the vanishing of the first Chern class is needed to ensure a well-defined integer grading). Moreover, note that the first condition is harmless, as otherwise we would obtain infinitely many fixed points by compactness. Also note that we make no assumption on non-degeneracy of periodic orbits.

\medskip

\textbf{Open-string version.} We also have the following relative version, which generalises a result of the authors \cite{LM2}. Given a Hamiltonian map $f$ of a symplectic manifold $W$, and a Lagrangian $L\subset W$, a pair $(x,m)$ where $x \in L$, $m \in \mathbb{N}\backslash\{0\}$, and $f^m(x) \in L$ is called a \textit{chord} of order $m$ of $\tau$. We call \textit{minimal order} of $x$ the smallest such $m$. A chord is periodic of period $k$ if $f^k(x) = x$, and we define its \textit{minimal period} as the smallest such $k$. A periodic chord is a finite collection of \textit{sub-chords} $(x_0, m_0), \dots, (x_l, m_l) = (x_0, m_0)$, where each $(x_i, m_i)$ is obtained from the previous one by applying $f^{m_i -1}$.

\begin{theorem}[\textbf{Open strings: Long interior chords}]\label{thm:longintchords}

Let $f:(W,\omega)\rightarrow (W,\omega)$ be a $C^0$-Hamiltonian twist map on a Liouville domain. Let $L\subset (W,\lambda)$ be an exact, spin, Lagrangian with Legendrian boundary. Assume the following:
\begin{itemize}
\item \textbf{(Periodic chords)} There are finitely many periodic chords;
\item \textbf{(First Chern class)} $c_1(W)=0$ if $\dim W \geq 4$;
\item\textbf{(Wrapped Floer cohomology)} $HW^\bullet(L)$ is non-zero in infinitely many degrees.
\end{itemize}

\smallskip

Then $f$ admits infinitely many interior chords with respect to $L$, of arbitrary large order, and which are not sub-chords of any periodic chord.

\end{theorem}

Note that we make no assumption on the non-degeneracy of chords. Similarly to the closed version above, in the open case we can get rid of the assumptions on index-definiteness (imposed in \cite{LM2}), we can drop the assumption that the contact structure at the boundary should be globally trivial, as well as relax the regularity from $C^2$ to $C^1$. The proof of Theorem \ref{thm:longintchords} is completely analogous to that of Theorem \ref{thm:main_thm}, and so we will not include it, referring the reader to \cite{L} for full details.

\smallskip

\textbf{Acknowledgements.} For this work, the authors were supported by the Air Force Office of Scientific Research (AFOSR) under Award No.\ FA8655-24-1-7012, by the DFG under Germany's Excellence Strategy EXC 2181/1 - 390900948 (the Heidelberg STRUCTURES Excellence Cluster), and by the Sonderforschungsbereich TRR 191 Symplectic Structures in Geometry, Algebra and Dynamics, funded by the DFG (Projektnummer 281071066). \

\section{Hamiltonian twist maps on degenerate Liouville domains}\label{sec:degenerate} In this section, we illustrate the trade-off alluded to in the introduction, and make the boundary degeneracy of the symplectic form explicit. We first introduce the notion of a degenerate Liouville domain, which roughly speaking is a Liouville domain away from the boundary, but whose symplectic form degenerates along the latter.

\begin{definition}\label{def:deg_Liouville}
A \emph{degenerate} Liouville domain is a triplet $(W,\lambda,\alpha)$ where:
\begin{itemize}
    \item $W$ is a smooth manifold with boundary $B=\partial W$;
    \item $\alpha$ is a contact form on $B$; and
    \item $\lambda$ is a $1$-form such that $\omega=\mathrm{d}\lambda$ is symplectic in the interior int$(W)$, but degenerates at the boundary $B$ along a normal direction.
\end{itemize}
The last condition means that there is a collar neighbourhood $C=[0,1]\times B\subset W$ of the boundary, with collar coordinate $r\in [0,1]$ such that $B=\{r=1\}$, along which $\lambda=A(r)\alpha$, where $A$ is a smooth function satisfying
\begin{itemize}
    \item $\partial_rA>0$ for $r<1$;
    \item $\partial_rA\vert_{r=1}\equiv 0$;
    \item $A\vert_{r=1}\equiv 1$.
\end{itemize}
The $1$-form $\lambda$ is a \emph{degenerate} Liouville form.
\end{definition}

\begin{remark}Note that the notion of \emph{ideal} Liouville domains due to Giroux is different from the above. Indeed, in the ideal case, the Liouville form has a pole at the boundary, which is not the case here. Moreover, while degenerate Liouville domains have finite symplectic volume, ideal ones have infinite volume. 
\end{remark}

Sources of degenerate Liouville domains are pages of open book decompositions, and domains of billiard maps, as explained in what follows.

\begin{example}\label{examples}
\begin{itemize} $\;$
    \item[(1)] \textbf{(Open books)} Let $(M,\xi=\ker \alpha)$ be a contact manifold where $\alpha$ is a contact form adapted to an open book decomposition $(B,\theta)$, $\theta: M \backslash B\rightarrow S^1$. This means that $\alpha_B:=\alpha\vert_B$ is a positive contact form, and $\mathrm{d}\alpha$ is symplectic in the fibres of $\theta$. We assume that the open book is of abstract type $M=\mathbf{OB}(W,\phi)$ where $W$ is the Weinstein page (the closure of the typical fibre of $\theta$), and $\phi: W\rightarrow W$ is the monodromy of the fibration $\theta$, with $\phi\vert_B=\text{id}$. 
    
    Let $\omega=\mathrm{d}\lambda$ with $\lambda=\alpha\vert_{\mbox{int}(W)}$. Then $\omega\vert_B$ is degenerate, as the Reeb vector field of $\alpha$ is tangent to $B$. Moreover, \cite[Prop.\ 3.1]{DGZ} provides a neighbourhood $B\times \mathbb{D}^2\subset M$ of the binding in which $\alpha=A(\alpha_B+s^2 \mathrm{d}\theta),$ where $(s,\theta)$ are polar coordinates, for a smooth positive function $A$ satisfying $A\equiv 1$ along $s=0$ and $\partial_sA<0$ for $s>0$. The contact condition also implies that $A=A(s)$ only depends on $s$, and the fact that $R_\alpha$ is tangent to $B$ implies $\partial_sA\vert_{s=0}=0$. Then $\lambda=A(s)\alpha_B$, and by changing coordinates to $r=1-s$, we see explicitly that $(W,\omega)$ is a degenerate Liouville domain. 

    There is a natural symplectic map $f: (\mbox{int}(W),\omega)\rightarrow (\mbox{int}(W),\omega)$, given by the Poincar\'e return map of $\alpha$. In general this map may not extend to the boundary, but it does whenever the Hessian of $A$ is negative-definite; see \cite{MvK}.  
    \smallskip

    \item[(2)] \textbf{(The Circular Restricted Three-Body Problem)} In \cite{AFvKP, CJK}, it was shown that the regularised flow of the spatial Circular Restricted Three-Body Problem (CR3BP) could be viewed as a Reeb flow on a contact $5$-fold, for sub-critical energies and near the primaries. In \cite{MvK2} it was shown that this $5$-fold admits an open book decomposition adapted to the flow. Therefore, we have a family of global hypersurfaces of section for the CR3BP, which by the previous point are all degenerate Liouville domains; with a return map which extends smoothly to the boundary by \cite{MvK2}.

    \smallskip
    \item[(3)] \textbf{(Billiards)} Let $D\subset \mathbb R^2$, $D\cong \mathbb D^2$, be a planar and strictly convex billiard table. The associated Birkhoff annulus $W\cong [-1,1]\times S^1$ is the set of all unit vectors pointing inwards along $\partial D$, which can be parametrised by two variables $\theta \in [0,1]$ and $\varphi \in S^1$. $W$ admits the two-form $\omega=\mathrm{d}(-\cos(\theta)\mathrm{d}\varphi)=\sin(\theta)\mathrm{d}\theta\wedge \mathrm{d}\varphi$, so that $(W,\omega)$ is a degenerate Liouville domain. The associated symplectic billiard map $f:(\mbox{int}(W),\omega)\rightarrow (\mbox{int}
(W),\omega)$, mapping a vector to the next vector along a billiard trajectory, extends smoothly to the boundary as the identity. In a different description where the $2$-form is non-degenerate at the boundary, the billiard map is only continuous at the boundary.
\end{itemize}   
\end{example}

\textbf{Degeneration and non-degeneration.} A degenerate Liouville domain $(W,\omega=\mathrm{d}\lambda)$ can be turned into a Liouville domain $(W,\omega_Q=\mathrm{d}\lambda_Q)$, which is unique up to Liouville isotopy. We will call $(W,\omega_Q)$, the \emph{non-degeneration} of $(W,\omega)$. Similarly, a Liouville domain has an associated \emph{degeneration}. More precisely:

\begin{lemma}\textbf{(Degeneration and non-degeneration)}\label{lemma:degeneration} We have the following.
\begin{itemize}
    \item[(1)]\textbf{(From degenerate to non-degenerate)} Let $(W,\lambda,\alpha)$ be a degenerate Liouville domain. Then there exists a boundary-preserving homeomorphism $Q:W\rightarrow W$ which is smooth in $\mbox{int}(W)$ and only continuous at $B=\partial W$, such that:

    \begin{itemize}
        \item[$\bullet$]  $(W,\lambda_Q=Q^*\lambda)$ is a Liouville domain with strict contact-type boundary $(B,\alpha)$. This Liouville domain is unique up to Liouville isotopy. 
        \item[$\bullet$]   Given a smooth map $f:(W,\lambda,\alpha)\rightarrow (W,\lambda, \alpha)$, then $f_Q=Q^{-1} \circ f \circ Q$ is smooth in the interior, where it preserves $\omega_Q$, but only continuous at $B$.
    \end{itemize}
    
\item[(2)] \textbf{(From non-degenerate to degenerate)} Let $(W,\lambda)$ be a Liouville domain with strict contact boundary $(B,\alpha)$. Then, there exists a boundary-preserving homeomorphism $S=Q^{-1}:W\rightarrow W$ which is smooth, and such that $(W,\lambda_S=S^*\lambda, \alpha)$ is a degenerate Liouville domain. 
   
\end{itemize}

\[\begin{tikzcd}
	{(W, \lambda_Q = Q^\star\lambda_S)} &&& {(W, \lambda_S = S^\star\lambda_Q,\alpha)}
	\arrow["Q", shift left, from=1-1, to=1-4]
	\arrow["S = Q^{-1}", shift left, from=1-4, to=1-1]
	\arrow["f_Q = Q^{-1}\circ f \circ Q", loop left, from=1-1, to=1-1]
    \arrow["f", loop right, from=1-4, to=1-4]
\end{tikzcd}
\begin{tikzpicture}[overlay, remember picture]
	\node at (-9.4,-0.5) {\footnotesize Non-degenerate};
    \node at (-3.6,-0.5) {\footnotesize Degenerate};
\end{tikzpicture}
\]

\medskip

\end{lemma}

Here, $S$ is called a \emph{squaring} map, and $Q$ is called a \emph{square root} map. Both are reparametrisations in the direction normal to the boundary, i.e.\ the Liouville direction. It follows that if $f$ is generated by $H_t$ on the interior, $f_Q$ is generated by $H_t^Q=H_t\circ Q$ on the interior, but the isotopy is not-well defined at the boundary. See Definition \ref{def:C0HamtwistmapsSecondDef} below.

\medskip

To put it simply, the trade-off expressed in Lemma \ref{lemma:degeneration} is that in general we can either choose coordinates in which the $2$-form is symplectic but the map only continuous at the boundary, or alternatively where the $2$-form degenerates at the boundary but the map is smooth.

\begin{proof}[Proof of Lemma \ref{lemma:degeneration}] Take a collar $[0,1]\times B$ where $\lambda=A(s)\alpha$ is as in Definition \ref{def:deg_Liouville}. We consider a map $Q: W \rightarrow W$ which is the identity away from $[0,1]\times B\subset W$, and on $(0,1]\times B$ is of the form $Q(s,b)=(\varphi(s),b)$ where $\varphi$ solves the ODE
$$
\varphi^\prime(s)=-\frac{1}{\partial_sA(\varphi(s))}>0, \mbox{ for } s>0,
$$
which has a smooth solution for $s>0$. Integrating this equation, we see that
$$
A(\varphi(s))=1-s,
$$
for $s>0$, by choosing the integration constant to be $1$. By continuity, we see that this equation has a unique solution satisfying $\varphi(0)=0$, which is only continuous at $s=0$. Therefore $Q$ is a diffeomorphism in the interior of $W$, but it extends only continuously to the boundary. 

If we define $\lambda_Q=Q^*\lambda, \omega_Q=Q^*\omega=\mathrm{d}\lambda_Q$, then $\omega_Q$ is a symplectic form satisfying $\omega_Q=\mathrm{d}((1-s)\alpha)$. If we make the change of coordinates $r=1-s$, we see that $\omega_Q=\mathrm{d}(r\alpha)$ near $B=\{r=1\}=\{s=0\}$, and therefore $(W,\omega_Q)$ is a Liouville domain with strict contact-type boundary $(B,\alpha)$. Note that the maps $Q,Q^{-1}$ are given in the $r$-coordinate as
$$
Q(r,b)=\left(F(r),b\right),\;Q^{-1}(r,b)=(A(r),b),
$$
where $F(r)=1-\varphi(1-r)$. 

The only choice we made above consisted of the coordinates in which $\lambda=A(s)\alpha$. Any other choice of such coordinates differs by a reparametrisation in the $s$-direction, which doesn't change the Liouville isotopy class of the resulting domain. The statement about self-maps is obvious.

For the converse, write $\lambda=(1-s)\alpha$ where $s=1-r$ on $[0,1]\times B$. Let $S:W\rightarrow W$ which is the identity away from $[0,1]\times B\subset W$, and on $(0,1]\times B$ is of the form $S(s,b)=(s^2,b)$. Then $\lambda_S=S^*\lambda=(1-s^2)\alpha$ defines a degenerate Liouville form. Note that to remove it, the corresponding map $Q$ is $Q(s)=\sqrt{s}$. This finishes the proof.
\end{proof}

Note that from the proof of the Lemma, we see that a squaring map is of the form $Q=(F,\text{id})$ near $B$, where $F$ is strictly positive and equal to $1$ at the boundary. Moreover, the derivative $F'$ is strictly positive in the interior, and it explodes at the boundary.

\medskip

\textbf{Weakened Twist Condition.} The following is the degenerate version of a $C^0$-Hamiltonian twist map, obtained by taking Definition \ref{def:C0Hamtwistmaps} through the degeneration process.

\begin{defn}\textbf{(Hamiltonian twist maps on degenerate Liouville domains)}\label{def:C0HamtwistmapsSecondDef} Let $f:(W,\omega)\rightarrow(W,\omega)$ be a smooth map on a degenerate Liouville domain. We say that $f$ is a $C^0$-\emph{Hamiltonian twist map}, if the following hold:
\begin{itemize}
    \item\textbf{(Hamiltonian)} $f\vert_{\mbox{int}(W)}=\phi_H^1$ is generated by a $C^1$ Hamiltonian $H_t:\mbox{int}(W)\rightarrow \mathbb{R}$;
    \item \textbf{(Extension)} The Hamiltonian $H_t$, together with its Hamiltonian vector field, admit a $C^1$ extension to the boundary; and
    \item\textbf{(Weakened Twist Condition)} At the boundary $B$, the generating Hamiltonian vector field satisfies $h_t:=\alpha(X_{H_t})>0$.
\end{itemize}
 
We say that the generating Hamiltonian isotopy $H_t$ is \emph{strictly twisting} or \emph{positively wrapping}.
\end{defn}

We now consider an explicit local model near the boundary.

\begin{remark}\textbf{(Local models)} If $(W,\omega)$ is a Liouville domain with $\omega=\mathrm{d}(r\alpha)$ near the boundary $B$, and $E_t:(W,\omega)\rightarrow \mathbb{R}$ is an arbitrary $C^1$ Hamiltonian, its Hamiltonian vector field near $B$ is given by
\begin{equation}\label{eq:general}
X_{E_t}=(\partial_rE_t)R_\alpha +\frac{1}{r}\left(X_{E_t}^\xi-\mathrm{d}E_t(R_\alpha)V\right),
\end{equation}
where $R_\alpha$ is the Reeb vector field of $\alpha$, $V=r\partial_r$ is the Liouville vector field in $W$, and $X_{E_t}^\xi$ is defined implicitly via $i_{X_{E_t}^\xi}\mathrm{d}\alpha=-\mathrm{d}E_t\vert_{\xi}$. 

Now, say that instead, $E_t$ is a $C^1$ Hamiltonian on the degeneration $(W, \omega_S)$ of $(W,\omega)$ (given by the square root map $Q$, of the form $Q = (F, \text{id})$ near $B$). Consider the composition $H_t := E_t \circ Q$. In other words:

\[\begin{tikzcd}
	{(W, \omega)} &&& {(W, \omega_S)} &&& {\mathbb{R}}
	\arrow["Q", from=1-1, to=1-4]
	\arrow["H_t"', curve={height=18pt}, from=1-1, to=1-7]
	\arrow["E_t", from=1-4, to=1-7]
\end{tikzcd}\]

so that $H_t$ should be understood as the non-degeneration of $E_t$. Noting that the derivatives of $E_t$ along $B$ coincide with those of $H_t$, we have
\begin{equation}\label{eq:XQ}
X_{H_t}=\left[((\partial_rE_t)\circ Q)\cdot F^\prime(r)\right]R_\alpha +\frac{1}{F(r)}\left[\left(X_{E_t}^\xi-\mathrm{d}E_t(R_\alpha)V\right)\circ Q\right].
\end{equation}
This is clearly not defined at $r=1$ since $F^\prime(r)$ gives a pole for the first summand at $r=1$. The second summand remains bounded at $r=1$, and the first summand gives the infinite wrapping whenever $\partial_rE_t$ is strictly positive.
\end{remark}

\section{Extension}\label{sec:ExtensionOfH}

In what follows, following \cite{MvK}, we will deal with the fact that the Hamiltonians in the definitions above are not necessarily admissible for symplectic cohomology, by constructing suitable extensions to the completion of the Liouville domain. 

Let $(W,\omega=\mathrm{d}\lambda)$ be a Liouville domain with completion $(\widehat W=W\cup [1,+\infty)\times \partial W,\widehat \omega=\mathrm{d}\widehat\lambda)$, $r$ be the coordinate in the cylindrical end, $B=\partial W$, $\alpha=\lambda\vert_B$, and $H=H_t$ be an arbitrary $C^1$ Hamiltonian. By shifting by a constant, which does not affect the dynamics, we may assume that $H_t$ is positive. The symplectic form on the cylindrical end is $\omega=\mathrm{d}(r\alpha)$, so we get $$X_{H_t}\vert_B=h_tR_\alpha+Z,$$ with $\partial_rH_t\vert_{r=1}=h_t$ (which we assume arbitrary for the moment), and $Z$ linearly independent of $R_\alpha$. The family of Hamiltonians $H_t$ is \textit{a priori} unsuitable to compute symplectic homology. To deal with this, following \cite{MvK}, we will construct an extension $\widehat H$ to the cylindrical end of $\widehat W$ that is linear at infinity. 

We construction this extension as follows.
Let $$H_0(b,t):=H_t\vert_{r=1},\;H_1(b,t):=\frac{\partial H_t}{\partial r}\Big\vert_{r=1}=h_t,$$ and let $R$ be the remainder:
$$
R=(H-(H_0+(r-1)H_1)) \frac{2}{(r-1)^2},
$$
which is defined near $B$. By construction $R$ is a smooth function on a half-space.
The functions $H_0$ and $H_1$ are $r$-independent, so admit obvious extensions to $r>1$, but the function $R$ is $r$-dependent, so we use \cite{S} to extend it to $r>1$. Call this extension $\overline R$. 
Now choose $\delta_1>\delta_0 >0$ and choose a decreasing cutoff function $\rho$ with $\rho|_{[1,1+\delta_0]}=1$ and $\rho(r)=0$ for $r>1+\delta_1$;
\begin{itemize}
\item put  $\widehat R(r,b,t) = \overline R(r,b,t) \cdot \rho(r)$;
\item put $\widehat H_j(r,b,t) =C_j$ for $r\geq 1+\delta_1$, where $C_j \geq \max H_j$, and $C_1>0$;
\item and put $\widehat H_j(r,b,t) = H_j(b,t) \cdot \rho(r) +(1-\rho(r)\, )C_j$, for $j=0,1$.
\end{itemize}
The extension is then defined as
\begin{equation}
\label{eq:Ht}
\widehat H:=\widehat H_0(r,b,t) +(r-1) \widehat H_1(r,b,t)+\frac{(r-1)^2}{2!}\widehat R(r,b,t)
.
\end{equation}

With these choices, we have that $\widehat H = C_1(r-1)+C_0$ for large $r$.
The extension $\widehat H$ is therefore linear at infinity, and by perturbing the slope $C_1$ we can assume that $C_1\notin \mathrm{spec}(\alpha)$. The same can be arranged for all iterates $\widehat H^{\# k}$ by possibly changing $C_1$. The resulting Hamiltonians are then all linear at infinity, but they may have $1$-periodic orbits that are degenerate. If all the degenerate $1$-periodic orbits are isolated, then we can still define the Floer cohomology $HF^\bullet(\widehat H^{\# k})$. Since this can be achieved by perturbation (see Lemma 4.6 in \cite{MvK}), we assume that this is the case. We then obtain, similarly as in \cite{MvK}, that:

\begin{proposition}\label{prop:limit} If $\widehat H^{\# k}$ are admissible extensions of $C^1$-Hamiltonian twist maps, then
$$
SH^\bullet(W)=\varinjlim_{k} HF^\bullet(\widehat H^{\# k}).
$$
    
\end{proposition}

\section{Action growth} As the extensions in the previous section potentially introduce undesired orbits along a collar neighbourhood of the boundary, we need a way of distinguishing trajectories in the interior $\text{int}(W)$ from trajectories on the cylindrical end $[1,+\infty)\times\partial W$. In \cite{MvK} and \cite{LM2}, this was done by using an \textit{index growth} argument (which we review and generalise in Appendix \ref{app:index}).

\indent This argument relied on showing that the Conley-Zehnder index of trajectories of length $T$ on the end $[1,+\infty)\times\partial W$ satisfied an estimate of the type: $$\mu_{CZ} > c\cdot T + d$$ where $c > 0, d \in \mathbb{R}$, and $\mu_{CZ}$ denotes the Conley-Zehnder index in some (global) trivialisation. This phenomenon allows to tell apart interior chords from the ones on the cylindrical end, as the indices of the latter grow to infinity under iteration and hence do not contribute to the symplectic cohomology in the limit. However proving it requires strong assumptions: that the Hamiltonian $H$ satisfies the twist condition as in \cite{MvK}, and that the boundary $(\partial W,\alpha:=\lambda|_{\partial W})$ is strongly index-definite (which further imposes that the contact structure should be trivial as a bundle).

\medskip

In this section, we show that the index growth argument can be replaced by an \textit{action growth} argument, which only assumes a quantitative version of the (weakened) twist condition. As we shall see, this setting still allows to prove a Poincaré-Birkhoff theorem (both for symplectic and wrapped Floer homology) and hence reproduce the work from \cite{MvK} and \cite{LM2}. Moreover, this quantitative condition will be for free in the cases of interest.

\medskip

Let $(W,\omega = \mathrm{d}\lambda)$ be a Liouville domain with boundary $B$, and $f : W \to W$ a $C^1$-Hamiltonian twist map. In what follows, we will further need to make the assumptions:

\medskip

\textbf{Assumptions.} \textbf{(Quantitative twist condition)}
    $f$ (or $H_t$) is said to satisfy the \emph{quantitative} twist condition if it can be generated by a $C^1$-Hamiltonian $H_t : W\to\mathbb{R}$ such that

    \smallskip

    \begin{enumerate}
        \item $H_t|_{B} > 0$,

        \item $\min\limits_{B} h_t > \max\limits_{B} H_t$.
    \end{enumerate}

    \smallskip
    Here, $$h_t = \partial_r H_t=\alpha(X_{H_t})=\langle X_{H_t}, R_\alpha\rangle,$$ with $r$ the Liouville coordinate near the boundary, so that $V = r\partial_r$ is the Liouville vector near $B$. Condition (2) is a quantitative version of the weakened twist condition defined in the Introduction.

\smallskip

Note that the quantitative twist condition cleary implies the weakened twist condition. As it turns out, the quantitative twist condition will be for free in the case where the $C^1$-Hamiltonian twist map comes from smoothing a $C^0$-Hamiltonian twist map, as in Section \ref{sec:smoothing} (which is the only case of interest).

\begin{proposition}[Action growth]\label{prop:actiongrowth}
    Let $(W, \omega = \mathrm{d}\lambda)$ be a Liouville domain with boundary $B$, and $H_t : W\to \mathbb{R}$ a Hamiltonian satisfying the quantitative twist condition. Then we can construct an extension $\widehat{H}$ of $H$ to the completion $\widehat W$ such that $\widehat H$ is linear at infinity, i.e 

    \vspace{-0.8em}

    \begin{equation*}
        \widehat{H} = ar - \varepsilon
    \end{equation*}

    \smallskip
    
    \hspace{-1.1em} and there exist constants $c > 0, d \in \mathbb{R}$ such that for every trajectory $x : [0,T] \to [1,+\infty)\times B$ we have

    \vspace{-0.6em}

    \begin{equation*}
        \mathcal{A}_{\widehat H}(x) < -c\cdot T + d.
    \end{equation*}
    Moreover, the constant $c$ grows to infinity if $\partial_r H$ grows to infinity in $C^0$-norm.
\end{proposition}

In particular, the action of the trajectories in the cylindrical end goes to $-\infty$ linearly with their length $T$.

\begin{rem}
    As we shall see in the proof, we have quite some control over the linear extension $\widehat H$ we construct. The only constraints on the constants $a$ and $\varepsilon$ are

    \vspace{-1em}

\begin{align}
    a &> \max\limits_{B} \partial_r H \label{expr:conditionona}\\
    0 < \varepsilon &< \min\limits_{B} \partial_r H - \max\limits_{B} H \label{expr:conditiononespilon}
\end{align}

\medskip

Note that the right-hand side of (\ref{expr:conditiononespilon}) is positive by the quantitative twist condition. In particular, we can construct our extension to be arbitrarily close to the Hamiltonian $r \mapsto ar$.
\end{rem}

\textit{Proof of Proposition \ref{prop:actiongrowth}.} We explained, in Section \ref{sec:ExtensionOfH}, how to construct a linear extension $\widehat H$ of $H$ to the Liouville completion $\widehat W = W \cup_{B} [1,+\infty)\times B$, and we refer to that section for notation. Note that by construction $\widehat H$ admits no $1$-periodic orbits on $[1+\delta_1, \infty)\times B$, so that in particular it suffices to prove our desired action growth property on the neighbourhood $[1,1+\delta_1)\times B$.

\bigskip

\textbf{Step 1 (Approximating $\mathcal{A}_{\widehat H}$).} Recall that the action of a trajectory $x : [0,T] \to (\widehat W,\widehat \lambda)$ of $\widehat H$ is given by $$\mathcal{A}_{\widehat H}(x) = -\int_0^T x^\star\hat{\lambda} + \int_0^T \widehat H\big(x(t)\big)\mathrm{d}t.$$

\begin{rem}
    If we were in the presence of an exact Lagrangian, and were working with wrapped Floer cohomology, then the action functional would be the same up to the addition of a bounded term. In particular, our proof for symplectic cohomology will also carry out to the wrapped case.
\end{rem}

Since $x$ is a Hamiltonian trajectory, using Equation (\ref{eq:Ht}) we have

\vspace{-1em}

\begin{align*}
    x^\star\hat{\lambda} &= \hat{\lambda}(\dot{x}(t))\mathrm{d}t \\
    &= r\alpha(X_{\widehat H})\mathrm{d}t \\
    &= r\partial_r \widehat H \, \mathrm{d}t \\
    &= r\left(\partial_r \widehat{H}_0 + \widehat{H}_1 + \dots\right)
\end{align*}

where, when computing $\partial_r \widehat H$, we ignore all terms of order $\geq 1$ in $r-1$. This is because, when we chose our cut-off function $[1,1+\delta_1]\to[0,1]$, we had the freedom to choose $\delta_1 > 0$ arbitrarily small. In particular, we can ensure that the neighbourhood $[1,1+\delta_1)\times B$ is arbitrarily thin around $B$, so that we have $r \simeq 1$ everywhere. We also ignore the terms of order $\geq 1$ in the expression for $\widehat H$, so that $$\widehat H \simeq \widehat H_0$$ In particular, the action functional becomes:

\vspace{-1em}

\begin{align*}
    \mathcal{A}_{\widehat H}(x) &= -\int_0^T r\partial_r \widehat H\, \mathrm{d}t + \int_0^T \widehat H\big(x(t)\big)\mathrm{d}t \\ 
    &\simeq -\int_0^T r\big(\partial_r \widehat{H}_0 + \widehat H_1\big)\mathrm{d}t + \int_0^T \widehat{H}_0 \, \mathrm{d}t \\
    &\simeq \int_0^T \big(\widehat{H}_0 - \widehat{H}_1\big)\mathrm{d}t -\int_0^T \partial_r \widehat{H}_0 \, \mathrm{d}t
\end{align*}

\medskip

We can further simplify. Indeed:

\vspace{-1em}

\begin{align*}
    \widehat{H}_0-\widehat{H}_1 &= \rho(r)\Big((H_0-H_1) - (C_0-C_1)\Big) + C_0-C_1 \\
    \partial_r \widehat{H}_0 &= \rho'(r) \big(H_0-C_0\big)
\end{align*}

so that:

\vspace{-0.6em}

\begin{equation*}
    \mathcal{A}_{\widehat{H}}(x) \simeq \int_0^T \rho(r) \Big((H_0 - H_1) - (C_0 - C_1)\Big) \mathrm{d}t + \int_0^T (C_0-C_1)\mathrm{d}t - \int_0^T \rho'(r) \big(H_0-C_0\big) \,\mathrm{d}t.
\end{equation*}

\bigskip

\textbf{Step 2 (Estimating $\mathcal{A}_{\widehat{H}}$ in different zones).} The trick now is, instead of attacking the above expression heads-on, to first decompose $[1,1+\delta_1)\times\partial W$ into three zones, and approximate $\mathcal{A}_{\widehat{H}}$ in each of these. Recall that our cut-off function $\rho : [1,1+\delta_1]\to[0,1]$ is chosen to be decreasing, and flat near $1$ and $1 + \delta_1$, so that it looks like: 

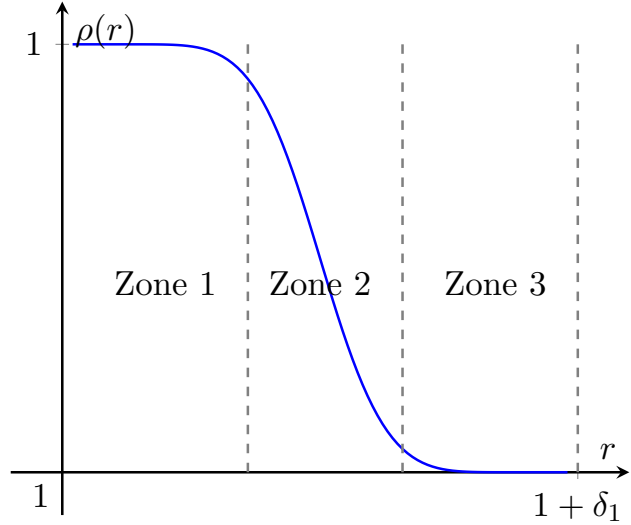
\[
    \begin{adjustbox}{scale=1.2}
   
    \begin{tikzpicture}
    \begin{axis}[
        domain=0.01:0.49, 
        samples=100,      
        axis lines=middle, 
        xlabel={$r$},     
        ylabel={$\rho(r)$},  
        ymin=0, ymax=1,   
        xmin=0, xmax=0.5, 
        grid=none,        
        thick,            
        xtick={0, 0.5},   
        ytick={0, 1},     
        xticklabels={0, $1 + \delta_1$}, 
        yticklabels={0, 1},        
        enlargelimits=true 
    ]
    \addplot[
        blue, 
        thick, 
    ]
    {(exp(1/x))/(exp(1/x) + exp(1/(-x + 0.5)))}; 
    
    \draw[dashed, gray] (axis cs:0.18,0) -- (axis cs:0.18,1);
    \draw[dashed, gray] (axis cs:0.33,0) -- (axis cs:0.33,1);
    \draw[dashed, gray] (axis cs:0.5,0) -- (axis cs:0.5,1);
    
    \node at (axis cs:0.1,0.5) [anchor=north] {Zone 1};
    \node at (axis cs:0.25,0.5) [anchor=north] {Zone 2};
    \node at (axis cs:0.42,0.5) [anchor=north] {Zone 3};
    
    \node at (axis cs:0,0) [anchor=north east] {1}; 
    \end{axis}
\end{tikzpicture}
\end{adjustbox}
\]

\medskip

In particular, we can decompose the interval $[1,1+\delta_1]$ into three zones, such that:

\begin{itemize}
    \item In Zone 1, $\rho(r) \simeq 1$, $\rho'(r) \simeq 0$.

    \item In Zone 2, $\rho$ decreases rapidly, i.e $\exists\, \tilde{\delta} > 0$ such that $\rho'(r) < -\tilde{\delta}$.

    \item In Zone 3, $\rho(r) \simeq 0 \simeq \rho'(r)$.
\end{itemize}

\medskip

We can then estimate $\mathcal{A}_{\widehat H}$ on each of these zones (which we choose to be open, and to slightly overlap, so that they cover the whole interval).

\medskip

\textbf{In Zone 1.} We get: $$\mathcal{A}_{\widehat{H}}(x) \simeq \int_0^T (H_0-H_1)\mathrm{d}t$$ In particular, by the quantitative twist condition we have $\min H_1 > \max H_0$, so that there must exist a constant $c_1 > 0$ such that 

\vspace{-0.6em}

\begin{equation}
    \mathcal{A}_{\widehat{H}}(x) \lesssim -c_1\cdot T.
\end{equation}

\medskip

\textbf{In Zone 3.} We get: $$\mathcal{A}_{\widehat{H}}(x) \simeq \int_0^T \big(C_0-C_1\big)\mathrm{d}t.$$ Therefore, as long as we choose our constants $C_0$ and $C_1$ so that $C_1 > C_0$ (which is possible, since the only constraint is that $C_i \geq \max H_i$), there exists a constant $c_3$ such that 

\vspace{-0.6em}

\begin{equation}
    \mathcal{A}_{\widehat H}(x) \lesssim -c_3\cdot T.
\end{equation}

\bigskip

\textbf{In Zone 2.} Zone 2 is the trickiest because none of the terms in $\mathcal{A}_{\widehat H}$ vanish. However, we can estimate each of them independently.

\begin{itemize}
    \item The last term $$-\int_0^T \rho'(r)\big(H_0-C_0\big) \mathrm{d}t$$ can be bounded by using the fact that $\rho'(r) < -\tilde{\delta}$ on Zone 2. And since we have control over the constant $C_0$, we can impose $C_0 > \max H_0$, so that there exists $c_2'' > 0$ s.t 

    \vspace{-0.6em}

    \begin{equation*}
        -\int_0^T \rho'(r)\big(H_0-C_0\big) \mathrm{d}t < -c_2''\cdot T. 
    \end{equation*}

    \item The second term is given by the integral of $C_0-C_1$, which we already showed was bounded by $-c_3\cdot T$.

    \medskip
    
    \item The last term is

    \begin{equation*}
        \int_0^T \rho(r) \Big((H_0-H_1) - (C_0-C_1)\Big) \mathrm{d}t
    \end{equation*}

    \smallskip

    In Zone $2$, the function $\rho$ has a minimum $\tilde{\varepsilon} > 0$, and if we choose the constants $C_0$ and $C_1$ such that

    \vspace{-0.6em}

    \begin{equation}\label{expr:conditionsOnConstants}
        C_0 < C_1 < C_0 + \min(H_1 - H_0)
    \end{equation}
    
    then we have ensured that there exists a $c_2' > 0$ such that 

    \vspace{-0.6em}

    \begin{equation*}
         \int_0^T \rho(r) \Big((H_0-H_1) - (C_0-C_1)\Big) \mathrm{d}t < -c_2'\cdot T.
    \end{equation*}
    \end{itemize} 

    \medskip

    \begin{rem}
        Note that, since $\min H_1 > \max H_0$ by the quantitative twist condition, we have $$\min H_1 - \max H_0 \leq \min(H_1 - H_0),$$ so that we can rewrite the condition (\ref{expr:conditionsOnConstants}) as 

        \vspace{-0.6em}

        \begin{equation*}
            C_1 = C_0 + \varepsilon \hspace{0.5em} \text{ for } 0 < \varepsilon < \min H_1 - \max H_0
        \end{equation*}
    \end{rem}

    \medskip

    In conclusion, in Zone 2, we have:

    \vspace{-0.8em}

\begin{align*}
    \mathcal{A}(x) &\simeq \underbrace{\int_0^T \rho(r) \Big((H_0 - H_1) - (C_0 - C_1)\Big) \mathrm{d}t}\limits_{\text{\normalsize $\lesssim -c_2'\cdot T$}} + \underbrace{\int_0^T (C_0-C_1)\mathrm{d}t}\limits_{\text{\normalsize $\lesssim -c_3\cdot T$}} - \underbrace{\int_0^T \rho'(r) \big(H_0-C_0\big) \,\mathrm{d}t}\limits_{\text{\normalsize $\lesssim -c_2''\cdot T$}} \\
    &\lesssim -c_2\cdot T,
\end{align*}

\medskip

by setting $c_2 = c_2' + c_2'' + c_3$.
\bigskip

\textbf{Step 3 (Patching the zones together).} Of course, there is no reason to assume that a trajectory $x : [0,T] \to [1,1+\delta_1)\times\partial W$ of the Hamiltonian flow will be constrained to one of the three zones. In general, it may travel freely between them. This can easily be fixed by taking a partition of the interval $[0,T]$, sub-dividing it into intervals along which $x$ is constrained to a specific zone (by which we mean, its $r$-coordinate belongs to one of Zone $1,2,3$). Then, by setting $c := \min(c_1,c_2,c_3)$, we get

\vspace{-0.6em}

\begin{equation*}
    \mathcal{A}_{\widehat{H}}(x) \lesssim -c\cdot T
\end{equation*}

The term $\lesssim$ means that this inequality is true up to a small error, so in other words, $\exists \, d \in \mathbb{R}$ s.t 

\vspace{-0.6em}

\begin{equation*}
    \boxed{\mathcal{A}_{\widehat{H}}(x) \leq -c\cdot T + d}
\end{equation*} 

 Note that in all Zones, we see that $c_i$ is large if $H_1$ (and hence $C_1$) is large, which happens if $\partial_rH$ is large. This proves the last statement, concluding the proof. \qed

\medskip

Recall that, in Step 2, we imposed the conditions $C_0 > \max H_0$, and $C_1 = C_0 + \varepsilon$ for some $\varepsilon < \min H_1 - \max H_0$; and we must also have $C_1 \geq \max H_1$ from the construction of our Hamiltonian. To summarise, our extension must be of the form: 

\vspace{-1em}

\begin{align*}
    \widehat{H} &= (C_0+\varepsilon)(r-1) + C_0 \\
    &= (C_0+\varepsilon)r - \varepsilon
\end{align*}

\smallskip

where $\varepsilon < \min H_1 - \max H_0$, and $C_0 \geq \max H_1 - \varepsilon$. Or more concisely, letting $a := C_0 + \varepsilon$, we obtain

\vspace{-0.4em}

\begin{equation*}
    \widehat{H} = ar - \varepsilon.
\end{equation*}

\begin{remark}
    The quantitative twist condition is key in proving the action growth of chords \textit{extremely close} to the boundary $B$ (in Zone 1). It is also important in Zone 2 (especially at the beginning), to make sure the first term satisfies action growth, but the more we progress along Zone 2, the more $\mathcal{A}_{\widehat H}$ will become dominated by the integral in $-\rho'(r)$, because $\rho$ will decrease rapidly. Finally, in Zone 3, we don't need the quantitative twist condition at all, and simply achieve action growth thanks to our choice of constants $C_0, C_1$.
\end{remark}

\section{Smoothing a $C^0$-Hamiltonian twist map}\label{sec:smoothing} The key idea for the proof of the main results is to approximate a given $C^0$-Hamiltonian twist map $f$ by a family of $C^1$-Hamiltonian twist maps $f_\epsilon$, where $f_\epsilon$ is generated by a $C^1$ Hamiltonian $H_\epsilon$ and converges to $f$ as $\epsilon \rightarrow 0$, in such a way that the slopes of the extensions $\widehat H_\epsilon$ grow to infinity. Then taking a limit, we compute the symplectic cohomology of $W$, and then the arguments of \cite{MvK} apply.

\smallskip

\begin{definition}\textbf{($C^1$-Hamiltonian twist map)} Let $f: (W,\omega)\rightarrow (W,\omega)$ be a map on a Liouville domain. We say that it is a \emph{$C^1$-Hamiltonian twist map}, if
\begin{itemize}
    \item\textbf{(Hamiltonian)} $f=\phi_H^1$ is generated by a $C^1$ Hamiltonian $H_t:W\rightarrow \mathbb R$;
    \item\textbf{(Weakened Twist Condition)} At the boundary $B$, the generating Hamiltonian vector field satisfies $h_t:=\alpha(X_{H_t})>0$.
\end{itemize}

We say that the isotopy $H_t$ is \emph{strictly wrapping} or \emph{positively wrapping}.

\end{definition}

\medskip

Let $f:(W,\omega)\rightarrow (W,\omega)$ be a $C^0$-Hamiltonian twist map on a Liouville domain, as in the statement of the main results, with infinitely wrapping generating Hamiltonian $H_t: W \rightarrow \mathbb R$. We will prove the following.

\begin{theorem}\label{thm:smoothing}
    For $\epsilon \geq 0$, there exists a family of $C^1$-Hamiltonian twist maps $f_\epsilon$ on a Liouville domain $(W,\omega_\epsilon)$, such that:
    \begin{itemize}
        \item $f_\epsilon$ is generated by a $C^1$ Hamiltonian $H_\epsilon = H_{t,\epsilon}$ which converges in $C^0$ to $H_t$ as $\epsilon \rightarrow 0$.
        \item Along $B$, the function $h_{t,\epsilon}=\alpha(X_{H_{t,\epsilon}})=\partial_r H_{t,\epsilon}$ diverges uniformly and monotonically as $\epsilon \rightarrow 0$, but all derivatives of $H_{t,\epsilon}$ in directions tangent to $B$ remain uniformly bounded. 
    \item As $\epsilon \rightarrow 0$, $\omega_\epsilon$ converges to $\omega$ in $C^\infty$.
    \item The completion $\widehat \omega_\epsilon$ on $\widehat W$ is independent of $\epsilon>0$, i.e.\ it symplectomorphic to $\widehat \omega$.
    \item The slope of the extension $\widehat H_\epsilon$ on $\widehat W$ is bounded below by $C/\epsilon$ with $C>0$, and so monotonically diverges as $\epsilon \rightarrow 0$.
    \end{itemize}
\end{theorem}

A direct corollary of the above is the following.

\begin{corollary}\label{cor:limit} Let $H=H_t: W\rightarrow \mathbb R$ be an infinitely wrapping Hamiltonian on a Liouville domain $W$ generating a $C^0$-Hamiltonian twist map, and let $H_\epsilon=H_{t,\epsilon}$ be as in Theorem \ref{thm:smoothing}. Then:
\begin{itemize}
    \item We have $$
\lim_{\epsilon \rightarrow 0}HF(\widehat W,\widehat{H}_\epsilon)=SH(W).
$$
\item For $\epsilon$ sufficiently small, $H_\epsilon$ satisfies the quantitative twist condition.
\end{itemize}

\end{corollary}

\begin{proof}[Proof of Theorem \ref{thm:smoothing}] We view $(W,\omega=\omega_Q=Q^*\omega_S)$ as the non-degeneration of $(W, \omega_S)$. Recall that the degeneration is given by the square root map $Q: (W,\omega_Q)\rightarrow (W,\omega_S)$ as in Lemma \ref{lemma:degeneration}. Like in the lemma, assume $$Q(r,b)=(F(r),b)=(1-\varphi(1-r),b)$$ on a collar $(1-\epsilon,1]\times B$, where $B=\{r=1\}$. In particular, $F'(r)=\varphi'(1-r)$ has a pole at $r=1$. We write $H_t=E_t \circ Q=E_t^Q$ where $E_t$ is a $C^1$ Hamiltonian on the degeneration $(W,\omega_S)$. (See the local model from Section \ref{sec:degenerate}).

We let $g_\epsilon:[0,1]\rightarrow (0,\infty)$, for $\epsilon>0$, be a family of positive $C^\infty$ functions such that $g_\epsilon(s)=g_0(s)$ for $\epsilon\leq s\leq 1$, where $g_0=\varphi^\prime$, and such that $g_\epsilon(0)=1/\epsilon$. Then the $g_\epsilon$ are smooth truncations of $g_0$ near $s=0$, and $g_\epsilon\rightarrow g_\infty$ in the $C^0$ topology. We let $$\varphi_\epsilon(s)=\int_0^sg_\epsilon(x)\mathrm{d}x,$$ which are smooth, $\varphi_\epsilon(0)=0$ for all $\epsilon$, and $\varphi_\epsilon\rightarrow \varphi_0:=\varphi$ in $C^0([0,1])$. We consider a smooth map
$$
Q_\epsilon:W\rightarrow W
$$
given by the identity away from the collar $(1-\epsilon,1]\times B$, and which coincides with
$$
Q_\epsilon(s,b)=(\varphi_\epsilon(s),b)
$$
near $B$, so that $Q_\epsilon\rightarrow Q_0:=Q$ in $C^0$. Let $$F_\epsilon(r)=1-\varphi_\epsilon(1-r)$$ denote the change of coordinates from the $s$-coordinate to the $r$-coordinate corresponding to $Q_\epsilon$ (so that $Q_\epsilon(r,b) = (F_\epsilon(r),b)$), satisfying $F_\epsilon\vert_{r=1}\equiv 1$ for all $\epsilon$, and converging to $$F_0(r):=F(r)=1-\varphi(1-r).$$ Finally, we define $H_{t,\epsilon}:(W,\omega_\epsilon)\rightarrow \mathbb R$ by $H_{t,\epsilon}=E_t\circ Q_\epsilon$, so that $H_{t,\epsilon}$ coincides with $H_t$ away from the collar, converges to $H_t=E_t^Q$ in $C^0$, and
$$
H_{t,\epsilon}(r,b)=E_t(F_\epsilon(r),b).
$$ Note that their Hamiltonian vector fields, when computed with respect to the symplectic form $\omega_\epsilon=Q_\epsilon^* \omega_S$ (which looks like $\mathrm{d}(F_\epsilon(r)\alpha)$ near the boundary), are given near $B$ by
\begin{equation}\label{eq:truncation}
X_{H_{t,\epsilon}}=\left[((\partial_rE_t)\circ Q_\epsilon)\cdot g_\epsilon(1-r)\right]R_\alpha +\frac{1}{F_\epsilon(r)}\left[\left(X_{H_t}^\xi-\mathrm{d}H_t(R_\alpha)Y\right)\circ Q_\epsilon\right].
\end{equation}
The effect is that the first summand in this vector field no longer has a pole at $r=1$, and hence $X_{H_{t,\epsilon}}$ is indeed smooth. Moreover, note that as $\epsilon$ goes to zero, the vector field becomes more and more collinear with the Reeb vector field. Moreover, while the symplectic form a priori depends on $\epsilon$, since $F_\epsilon$ is positive, the effect is to change the contact form at the boundary, but we have that their completions $\widehat \omega_\epsilon$ on $\widehat W$ are all symplectomorphic. That is, their completion is an $\epsilon$-independent $2$-form $\widehat \omega$ on $\widehat W$, so that the Liouville structure is fixed.

By assumption, we have
$$
X_{H_t}=h_t R_\alpha + Z,
$$
where $h_t>0$, $h_t\vert_{r=1}\equiv \infty$, and $Z$ is linearly independent of $R_\alpha$. If we look at the truncated version, Equation (\ref{eq:truncation}) implies $X_{H_{t,\epsilon}}$ is of the form
$$
X_{H_{t,\epsilon}}=h_{t,\epsilon}R_\alpha + Z_\epsilon,
$$
where $h_{t,\epsilon}>0$, and $$h_{t,\epsilon}\vert_{r=1}=\frac{1}{\epsilon}\cdot (\partial_rH_t)\vert_{r=1}\rightarrow \infty$$ uniformly and monotonically as $\epsilon\rightarrow 0$. In other words, $H_{t,\epsilon}$ generates a $C^1$-Hamiltonian twist map $f_\epsilon$. This admits the following interpretation: if $H_t$ strictly twists at the boundary with respect to $\omega$, then $H_t^Q$ gives a strict ``infinite twist'' at the boundary when computed with respect to $\omega$, and $H_{t,\epsilon}$ is a truncation that strictly twists up to $1/\epsilon$ times the original twisting of $H_t$. Moreover, making $\epsilon$ smaller makes $h_{t,\epsilon}$ uniformly larger in $C^0$ norm, while keeping all $B$ derivatives of $H_{t,\epsilon}$ uniformly bounded in $C^0$ norm. 

Moreover, we immediately see that the extension $\widehat H_\epsilon$ of $H_\epsilon$ to $\widehat W$ is of the form
$$
\widehat H_\epsilon = C^\epsilon_1 (r-1) + C_0^\epsilon,
$$
with $C_1^\epsilon \geq \max_{B,t}(h_{t,\epsilon})\geq \frac{1}{\epsilon}C$, with $C:=\min_{B,t} \partial_r H_t>0$. This finishes the proof.
\end{proof}

\section{Proof of the main theorems}\label{sec:proofOfMainThms}

We now prove Theorem \ref{thm:main_thm}, namely that if $f : (W,\omega)\to(W,\omega)$ is a $C^0$-Hamiltonian twist map whose fixed points are isolated, and if the symplectic homology $SH(W)$ is non-zero in infinitely many degrees, then $f$ admits infinitely many interior periodic points of arbitarily large period. 

\smallskip

Conceptually, the argument is analogous to the one in the main theorem of \cite{MvK} (which in turn is adapted from Ginzburg's arguments in the proof of the Conley conjecture \cite{G10}), except we replace the index growth argument by our newly found \textit{action} growth argument, in order to distinguish interior orbits from the ones on the collar. 

\begin{rem}
    We do not include the proof of Theorem \ref{thm:longintchords}, the reason being that it is almost exactly the same as that of Theorem \ref{thm:main_thm}, simply replacing periodic orbits by Hamiltonian chords, and symplectic by wrapped Floer cohomology. For a full write-up of the proof of Theorem \ref{thm:longintchords}, see \cite{L}.
\end{rem}

\medskip

\textit{Proof of Theorem \ref{thm:main_thm}.} Let $H=H_t : W \to \mathbb{R}$ be the $C^0$-twist Hamiltonian generating the map $f$. By assumption, the fixed points of $f$ are isolated and hence there are finitely many, so let us denote them $\gamma_1, \dots, \gamma_k$. Assume for a contradiction that $f$ has finitely many interior periodic points $x_1,\dots,x_\ell$, with minimal periods $m_1, \dots, m_\ell$. For $\epsilon\geq 0$, denote by $H_\epsilon$ a smoothing of $H$ as in Section \ref{sec:smoothing}, and by $\widehat H_\epsilon$ the corresponding admissible extension. We assume that $\epsilon$ has been chosen small enough so that the $\gamma_j$ and the $f$-orbits of all $x_i$ lie away from the collar neighbourhood of the boundary where $H_\epsilon$ is actually $\epsilon$-dependent, so that they are still fixed points (resp.\ periodic points) of the $C^1$-Hamiltonian twist map $f_\epsilon$ generated by $H_\epsilon$. 

\smallskip

Choose a sequence of primes $\{p_i\}$ going to $+\infty$, and such that $p_i > \max_{k} m_k$. Let $\epsilon_i=1/p_i>0$. Using Proposition \ref{prop:limit} and Corollary \ref{cor:limit}, we obtain that

    \vspace{-0.6em}

    \begin{equation*}
        SH^{\bullet}(W) := \varinjlim\limits_{p_i} \, \, HF^{\bullet}(\widehat{H}_{\epsilon_i}^{\# p_i}), 
    \end{equation*}

    \noindent where $\widehat{H}_{\epsilon_i}^{\# p_i}$ is the $p_i$-th iterate of $\widehat H_{\epsilon_i}$. 

    \medskip

    Now pick $N > 2nk$, where $2n = \dim W$.  Recall that, by assumption, $SH^{\bullet}(W)$ is non-zero in infinitely many degrees, so that we can find $i_1, \dots, i_N$ such that $SH^{i_j}(W) \ne 0$. Therefore $$HF^{i_j}(\widehat{H}_{\epsilon_i}^{\# p_i}) \ne 0 \hspace{0.5em} \text{for } i \text{ large enough and for all } j=1,\dots,N.$$

By action growth (Proposition \ref{prop:actiongrowth}), we have, for any $1$-periodic orbit $x$ of $\widehat{H}_{\epsilon_i}^{\# p_i}$ on the collar $[1,+\infty)\times\partial W$, that

    \vspace{-0.6em}

    \begin{equation}\label{NegActionAnakin}
        \mathcal{A}_{\widehat H_{\epsilon_i}^{\# p_i}}(x) \leq -c_i\cdot p_i + d,
    \end{equation}
    with $c_i \rightarrow +\infty$ as $i\rightarrow +\infty$. Moreover, for every $i$ there exists a local-to-global spectral sequence whose $E_1$ page is $$E_1^{pq}(\widehat H_{\epsilon_i}^{\# p_i})=\bigoplus_{f_i(p-1)<\mathcal{A}_{\widehat H_{\epsilon_i}^{\# p_i}}(x)<f_i(p)} HF^{p+q}_{loc}(x,\widehat H_{\epsilon_i}^{\# p_i})$$ made of local Floer cohomologies of the (potentially degenerate) $1$-periodic orbits $x$ of $H_{\epsilon_i}^{\# p_i}$. Here, $f_i: \mathbb N_0 \rightarrow \mathbb R$ is a decreasing function (i.e.\ we order the isolated $1$-periodic orbits by decreasing action, so that large columns $p$ correspond to orbits with very negative action). Moreover, the spectral sequence converges to the Floer homology (see \cite{KvK} or \cite{LM2}):

    \vspace{-0.6em}

    \begin{equation*}
        E_1(\widehat H_{\epsilon_i}^{\# p_i}) \implies HF^{*}(\widehat{H}_{\epsilon_i}^{\# p_i}), \; i.e.\ HF^{i_j}(\widehat{H}_{\epsilon_i}^{\# p_i}) \cong \displaystyle\bigoplus\limits_{p+q = i_j} E_\infty^{p,q}.
    \end{equation*}
    
    Since $HF^{i_j}(\widehat{H}_{\epsilon_i}^{\# p_i}) \ne 0$ for all $j$, there must be non-zero elements on the diagonals $p + q = i_j$ of $E_\infty$. If $p_i \gg \max_j|i_j|$, then by (\ref{NegActionAnakin}), orbits of $\widehat{H}_{\epsilon_i}^{\# p_i}$ on $[1,+\infty)\times\partial W$ will have action escaping to $-\infty$, and hence appear on columns with $p \gg 1$. Therefore such orbits cannot contribute to the diagonals $p +q = i_j$, and \textit{a fortiori} towards $HF^{i_j}(\widehat{H}_{\epsilon_i}^{\# p_i})$, so that only interior orbits contribute.

   \smallskip

    \indent Now recall that $p_i$ is prime, so that the  $m_1,\dots,m_\ell$ cannot divide $p_i$, ensuring that none of the periodic points $x_1,\dots,x_\ell$ contribute to $HF^{i_j}(\widehat{H}_{\epsilon_i}^{\# p_i})$. Therefore, the orbits we have found contributing to $HF^{i_j}(\widehat{H}_{\epsilon_i}^{\# p_i})$ must necessarily be iterates $\gamma_j^{p_i}$ of one of the $1$-periodic orbits $\gamma_j$. However, a well-known fact about the mean index $\Delta(x)$ of a periodic orbit $x$ (see e.g.\ \cite{G10}) is that it always stays at a distance $\leq n=\dim(W)/2$ of the Conley-Zehnder index, and so

    \vspace{-0.6em}

    \begin{equation*}
        \text{supp} \,HF_\text{loc}^{\bullet}(\gamma_j^{p_i},\widehat H_{\epsilon_i}^{\# p_i}) \subset [\Delta(\gamma_j^{p_i}) - n, \Delta(\gamma_j^{p_i}) + n] \mbox{ for all } j,
    \end{equation*}

    \smallskip

   \noindent  where the \textit{support} denotes the range of degrees in which the local cohomology is non-vanishing. Hence, each orbit $\gamma_j^{p_i}$ can contribute to at most $2n$ degrees in cohomology, so that by counting all of them, we have covered $2nk$ degrees. However, we had found $N > 2nk$ values for the degree in which the cohomology is non-zero. This yields a contradiction. 

   We conclude that there must exist infinitely many interior periodic points, of arbitrarily large order, which finishes the proof. \qed 

\appendix

\section{Index growth}\label{app:index}

In this section, for completeness (although it is not needed in the main body of the text to prove the main theorems), we show how one can still obtain index growth of Hamiltonian trajectories of the extension of Section \ref{sec:ExtensionOfH}, under milder assumptions as those of \cite{MvK}, and in the case of a $C^2$-Hamiltonian twist map as defined in Section \ref{sec:degenerate}, with the only difference that $C^1$ regularity is replaced with $C^2$ regularity. 

We say that a strict contact manifold $(Y,\xi=\ker \alpha)$ is \emph{strongly index-definite} if $\xi$ admits a symplectic trivialisation $\epsilon$, and there are constants $c>0$ and $d\in \mathbb{R}$ such that for every Reeb arc $\gamma:[0,T]\rightarrow Y$ of Reeb action $T=\int_0^T \gamma^*\alpha$ we have
    $$
    \vert\mu_{RS}(\gamma;\epsilon)\vert\geq c T+d,
    $$
where $\mu_{RS}$ is the Robbin--Salamon index.

\medskip

We consider a Liouville domain $(W,\omega)$ with boundary $B$, and a $C^2$-Hamiltonian twist map with positively wrapping generating Hamiltonian $H_t$.

\medskip

\textbf{Assumptions.} Assume the following:
\begin{itemize}
    \item[\textbf{(A1)}] The linearised Reeb flow equation is strongly index definite along $B$. 
    \item[\textbf{(A2)}] $H_1=\partial_r H_t = h_t>0$ along $B$, i.e.\ the $C^1$-Hamiltonian twist condition.
    \item[\textbf{(A3)}] The $C^0$-norm of $h_t$ at $B$ can be chosen arbitrarily large in comparison to the $C^0$-norm of $H_t\vert_B$ and that of its derivatives up to order 2 with respect to directions in $B$.
\end{itemize}

The third condition, which was not in place in \cite{MvK}, is in fact for free, as we only care about the case when $H_t$ comes from smoothing a Hamiltonian generating a $C^0$-Hamiltonian twist map, and therefore $h_t$ can be chosen as large as desired. Note also that this condition implies the quantitative twist condition, as it implies $\min_B h_t \gg \max_B H_t.$

We now linearise the Hamiltonian vector field of the extension $\widehat H$ of an arbitrary Hamiltonian $H$, following \cite{MvK}. The Hamiltonian vector field of $\widehat{H}_t$ near $B$ is given by
\begin{equation}\label{X_Hcylend}\small
X_{\widehat H}=F\cdot R_\alpha +\frac{1}{r}\left(X_\xi
-G\cdot Y \right),
\end{equation}
where
$$
F=\partial_r \widehat H_0 +\widehat H_1+ (r-1)\partial_r \widehat H_1
+(r-1)\widehat R +\frac{(r-1)^2}{2} \partial_r \widehat R,
$$
$$
X_\xi = X^\xi_{\widehat H_0}+(r-1) X_{\widehat H_1}^\xi+\frac{(r-1)^2}{2}X_{\widehat R}^\xi,
$$
and
$$
G=\mathrm{d}\widehat H_0(R_\alpha)+(r-1)\mathrm{d}\widehat H_1(R_\alpha)+\frac{(r-1)^2}{2}\mathrm{d}\widehat R(R_\alpha).
$$
Here, $Y=r\partial_r$ is the Liouville vector field, and $X_{h}^\xi \in \xi$ is the $\xi$-component of the contact Hamiltonian vector field $X_{h}=hR_\alpha + X_{h}^\xi$ of a Hamiltonian $h: B\rightarrow \mathbb{R}$, defined implicitly by the equation $\mathrm{d}\alpha(X^\xi_{h},\cdot)=-\mathrm{d}h\vert_{\xi}$. Note that the only difference between the above expressions and those in Section 4.2 of \cite{MvK} is the presence of the first terms in $X_\xi$ and $G$ above, which vanish in \cite{MvK} due to the assumption that the twist condition as defined there is satisfied.

To compute the linearisation, we choose the Levi-Civita connection for the metric $1/r^2 \cdot \mathrm{d}r\otimes \mathrm{d}r +\alpha\otimes \alpha +\mathrm{d}\alpha(\cdot,J\cdot)$. This connection has the following properties:
\begin{itemize}
\item $\nabla Y=0$;
\item $\nabla_{R_\alpha} R_\alpha=0$ and $\nabla_Y R_\alpha=0$;
\item $\nabla_X R_\alpha \in \xi$ for all $X\in \xi$.
\end{itemize}
With respect to this connection, we have
\begin{equation}
\label{linearisedXH}
\nabla X_{\widehat H}= F \nabla R_\alpha + \mathrm{d}F \otimes R_\alpha + \frac{1}{r^2}\mathrm{d}r \otimes (X_\xi-GY)+\frac{1}{r}(\nabla X_\xi-\mathrm{d}G \otimes Y).
\end{equation}
We group terms, and write
$$
\nabla X_{\widehat H}= L_0 + L_1,
$$
where
$$
L_0=F \nabla R_\alpha + \mathrm{d}F \otimes R_\alpha +\frac{1}{r^2}\mathrm{d}r \otimes (X_\xi-GY)
-\frac{1}{r^2} \mathrm{d}G(Y) \mathrm{d}r \otimes Y
+\frac{1}{r^2} \mathrm{d}r \otimes \nabla_Y X_\xi,
$$
and
$$
L_1=\frac{1}{r}\left(
\nabla^\xi X_\xi +\alpha \otimes \nabla_{R_\alpha} X_\xi
-R_\alpha(G) \alpha \otimes Y-\mathrm{d}^\xi G \otimes Y
\right).
$$
Here, $\nabla^\xi= P_\xi \nabla\vert_\xi,$ where $P_\xi$ is the orthogonal projection to $\xi$, and $\mathrm{d}^\xi=\mathrm{d}\vert_\xi$. We will show that, with respect to the decomposition $T\widehat W=\xi \oplus \langle Y,R_\alpha\rangle$, we have the matrix representations:
\[
L_0=
\left(\begin{array}{@{}c|c@{}}
  F \cdot \nabla^\xi R_\alpha
  & \begin{matrix}
  U & 0 \\
  V & 0 \\
  \end{matrix} \\
\hline
  \begin{matrix}
  0 & 0\\
  W & Z
  \end{matrix} &
  \begin{matrix}
     a & 0 \\
     b & c
  \end{matrix}
\end{array}\right), 
\quad
L_1=\frac{1}{r}
\left(\begin{array}{@{}c|c@{}}
  \nabla^\xi X_\xi
  & \begin{matrix}
  0 & U' \\
  0 & V' \\
  \end{matrix} \\
\hline
  \begin{matrix}
  W' & Z'\\
  0 & 0
  \end{matrix} &
  \begin{matrix}
     0 & a' \\
     0 & 0
  \end{matrix}
\end{array}\right)
\]

Moreover, a careful analysis completely analogous as to what is done in \cite{MvK} shows that both $L_0$ and $L_1$ lie in $\mathfrak{sp}(n)$, so that we can further constrain the entries of $L_0$ to be of the form:
\[
L_0=
\left(\begin{array}{@{}c|c@{}}
  F \cdot \nabla^\xi R_\alpha
  & \begin{matrix}
  U & 0 \\
  V & 0 \\
  \end{matrix} \\
\hline
  \begin{matrix}
  0 & 0\\
  V & -U
  \end{matrix} &
  \begin{matrix}
     a & 0 \\
     b & -a
  \end{matrix}
\end{array}\right) \in \mathfrak{sp}(2n).
\]

Moreover, we can write
$$
X_\xi=X^\xi_{\widehat H_0} + (r-1)X_\xi',
$$
$$
G=\mathrm{d}\widehat H_0(R_\alpha) + (r-1) G',
$$
with
$$
X_\xi'= X^\xi_{\widehat H_1} + \frac{r-1}{2} X_{\widehat R}^\xi,
$$
and
$$
G'=\mathrm{d}\widehat H_1(R_\alpha)+\frac{r-1}{2}d\widehat R(R_\alpha).
$$
This induces a further decomposition for $L_1$, of the form
$$
L_1=\frac{1}{r}L_1'+\frac{(r-1)}{r}L_1'',
$$
where
$$
L_1'=
\left(\begin{array}{@{}c|c@{}}
  \nabla^\xi X^\xi_{\widehat H_0}
  & \begin{matrix}
  0 & U'' \\
  0 & V'' \\
  \end{matrix} \\
\hline
  \begin{matrix}
  W'' & Z''\\
  0 & 0
  \end{matrix} &
  \begin{matrix}
     0 & a'' \\
     0 & 0
  \end{matrix}
\end{array}\right),
$$
which involves first derivatives in $B$ directions of the first summands $X_{\widehat H_0}^\xi$ and $d\widehat H_0(R_\alpha)$ appearing in $X_\xi$ and $G$ respectively, and
$$
L_1''=
\left(\begin{array}{@{}c|c@{}}
  \nabla^\xi X'_\xi
  & \begin{matrix}
  0 & U''' \\
  0 & V''' \\
  \end{matrix} \\
\hline
  \begin{matrix}
  W''' & Z'''\\
  0 & 0
  \end{matrix} &
  \begin{matrix}
     0 & a''' \\
     0 & 0
  \end{matrix}
\end{array}\right)
$$
which involves first derivatives of the remaining summands $X_\xi'$ and $G'$. 

We will then prove the following.

\begin{lemma}
    Under Assumptions \textbf{(A1)-(A3)}, the linearised Hamiltonian flow equation $\dot\psi=\nabla_\psi X_{\widehat H}$ of the extension $\widehat H$ is strongly index definite along the cylindrical end $[1,+\infty)\times B$. 
\end{lemma}

\begin{proof} Exactly as in Lemma 4.2 in \cite{MvK}, using \textbf{(A2)}, for $\delta_1$ small enough, the function $F$ will be strictly positive. Combined with \textbf{(A1)}, this implies that the first block of $L_0$ is strongly index definite, and then Appendix E in \cite{MvK} implies that the ODE $\dot \psi=L_0 \psi$ is strongly index-definite. Making $\delta_1$ smaller if necessary, observing that this makes all terms with a $r-1$ factor small, we can make the matrix $\nabla X_{\widehat H}$ arbitrarily $C^0$ close to the matrix $L_0+L_1'$. But now, using \textbf{(A3)}, we see that we can achieve that the $C^0$-norm of $F$ dominates the entries of $L_1'$, as these consist of first derivatives with respect to $B$ directions of $X_{\widehat H_0}^\xi$ and $\mathrm{d}\widehat H_0(R_\alpha)$ (so, second derivatives of $\widehat H_0$ with respect to $B$ directions), which can be taken bounded and small in comparison to $F$ by the explicit way in which the extension is constructed. By continuity of index growth (see Lemma 2.2.9 in \cite{U99}), we conclude the proof.
\end{proof}

\section{Non-example} In this section, as a sanity check, we check that the Katok examples introduced in \cite{K73}, of non-reversible Finsler metrics on $S^n$ with only finitely many simple closed geodesics, do \emph{not} satisfy the weakened twist condition. We follow the exposition of \cite{MvK}, using Brieskorn manifolds.

Let
$$
\Sigma^{2n-1} := \left\{ (z_0,\ldots,z_n) \in \C^{n+1} ~\Bigg|f(z)=~\sum_j z_j^2 = 0 \right\}
\cap 
S^{2n+1},
$$
equipped with the contact form $\alpha =\frac{i}{2}\sum_j z_j \mathrm{d}\bar z_j -\bar z_j \mathrm{d}z_j$.

We describe the case of $n=2m+1$ odd, the case of even $n$ is obtained by simply dropping the first coordinate. We consider the following unitary coordinate transformation:
$$
w_0=z_0,
w_1=z_1,
w_{2j}=\frac{\sqrt 2}{2} (z_{2j}+i z_{2j+1}),
w_{2j+1}=\frac{i \sqrt 2}{2} (z_{2j}-i z_{2j+1})
\text{ for }j=1,\dots,m.
$$
In $w$-coordinates, $\alpha$ is given by
$$
\alpha=\frac{i}{2} \sum_j w_j \mathrm{d}\bar w_j -\bar w_j \mathrm{d}w_j,
$$ and $f$ is given by
$$
f(w)=w_0^2+w_1^2-2i\sum_{j=1}^{m}w_{2j}w_{2j+1}.
$$
For a tuple $\epsilon=(\epsilon_1,\dots, \epsilon_m) \in (-1,1)^m$, let 
$$
H_\epsilon(w)=\Vert w \Vert^2+\sum_{j} \epsilon_j( 
|w_{2j}|^2-|w_{2j+1}|^2
).
$$
For $\epsilon$ sufficiently small, this function is positive, so we can define a perturbed contact form by
$$
\alpha_\epsilon = H_\epsilon^{-1} \cdot \alpha.
$$ 
The Reeb flow is given by
$$
(w_0,\dots,w_n)\longmapsto
(e^{2\pi it}w_0,e^{2\pi it}w_1,e^{2\pi it(1+\epsilon_1)}w_2,e^{2\pi it(1-\epsilon_1)}w_3,\dots, e^{2\pi it(1+\epsilon_m)}w_{n-1},e^{2\pi it(1-\epsilon_m)}w_{n}).
$$
This flow has only $n+1$ periodic orbits if all $\epsilon_j$ are rationally independent. 

We now consider the open book fibration given by
$$
\Theta: \Sigma^{2n-1} \longrightarrow \C, (w_0,w_1,\ldots,w_n) \longmapsto w_0.
$$
The zero set of $\Theta$ defines the binding, the pages are the sets of the form $P_\theta=\{\arg \Theta=\theta\}$, $\theta \in S^1$. The return map for the page $P_0=\Theta^{-1}(\R_{>0})\cong \mathbb{D}^*S^{n-1}$ is
\[
\begin{split}
\Phi: P_0&\longrightarrow P_0, \\
p=(r_0,w_1,w_2,w_3,\ldots,w_{n-1},w_{n}) & \longmapsto 
(r_0, w_1,e^{2\pi i\epsilon_1}w_2,e^{-2\pi i \epsilon_1}w_3,
\ldots,\\
&\phantom{\longmapsto (~} e^{2\pi i\epsilon_{m}}w_{n-1}, e^{-2\pi i\epsilon_{m}}w_n).
\end{split}
\]
Here, $w_0=r_0 \in \mathbb{R}_{\geq 0}$. If all $\epsilon_j$ are irrational and rationally independent, this map has only two interior periodic points:
\[
\begin{split}
p_0=& \left(\frac{1}{\sqrt{2}},\frac{i}{\sqrt{2}},0,\dots,0\right)\\
q_0=&\left(\frac{1}{\sqrt{2}},-\frac{i}{\sqrt{2}},0, \dots, 0\right).\\
\end{split}
\]

It follows from our main theorem that \emph{any} Hamiltonian generating $\Phi$ necessarily fails to satisfy the weakened twist condition. We will check this explicitly for the concrete generating Hamiltonian found in \cite{MvK}.

The symplectic form on the interior of the page $P_0$ is the restriction of $\omega_\epsilon = \mathrm{d}\alpha_{\epsilon}$, and let
$$
H(w)= \Vert w \Vert^2,\quad \Delta_\epsilon(w) =\sum_{j} \epsilon_j( 
|w_{2j}|^2-|w_{2j+1}|^2),
$$
so $H_\epsilon=H+\Delta_\epsilon$. The return map $\Phi$ is generated by the $2\pi$-flow of the vector field
$$
X=i\sum_{j=1}^m\epsilon_j\left(
w_{2j}\frac{\partial}{\partial w_{2j}}-\overline w_{2j}\frac{\partial}{\partial \overline w_{2j}}
-
w_{2j+1}\frac{\partial}{\partial w_{2j+1}}+\overline w_{2j+1}\frac{\partial}{\partial \overline w_{2j+1}}
\right).
$$
This vector field is tangent to the page. Indeed, it has no $\frac{\partial}{\partial \omega_0}$ summand, and one can also explicitly check that $\mathrm{d}H(X)=\mathrm{d}f(X)=0$, so that it is tangent to the sphere and the the zero set of $f$, i.e.\ to $\Sigma$. Moreover, as in \cite{MvK}, we have
\[
\begin{split}
\iota_X(\mathrm{d} H_\epsilon^{-1}\wedge \alpha +H_\epsilon^{-1} \mathrm{d}\alpha )& =
-\alpha(X) \mathrm{d} H_\epsilon^{-1}+H_{\epsilon}^{-1} \iota_X \mathrm{d}\alpha\\
&=-\Delta_\epsilon \mathrm{d} H_\epsilon^{-1} -H_\epsilon^{-1} \mathrm{d} \Delta_\epsilon= -\mathrm{d} (H_\epsilon^{-1} \Delta_\epsilon).
\end{split}
\]
This means that $X=X_{K_\epsilon}$ is the Hamiltonian vector field of $K_\epsilon:=H_\epsilon^{-1} \Delta_\epsilon$.

To check that this Hamiltonian does not satisfy the weakened twist condition, we evaluate $\alpha_\epsilon\vert_{P_0}$ on it, and check that the resulting function fails to be positive along the binding $w_0=0$. The result is
$$
\alpha_\epsilon(X)=H_\epsilon^{-1}(w)\cdot\alpha(X)=H_\epsilon^{-1}(w)\sum_{j=1}^m \epsilon_j\left(\vert w_{2j}\vert^2 -\vert w_{2j+1}\vert^2\right)=K_\epsilon(w).
$$
This expression has mixed sign along the binding. Indeed, we see that choosing $\vert w_2 \vert^2=\vert w_3 \vert^2=1$, and letting
$$
p_1=(0,0,w_2,0,0,\dots,0),\;p_2=(0,0,0,w_3,0,\dots,0),
$$
we have $p_j\in \partial P_0$ lies in the binding of $\Sigma$ for $j=1,2$, and $$K_\epsilon(p_1)=\frac{\epsilon_1}{1+\epsilon_1}>0,$$$$K_\epsilon(p_2)=-\frac{\epsilon_1}{1-\epsilon_1}<0.$$

This is indeed the expected situation.

\end{document}